%% file: main.tex
\documentclass{ifacconf}

\usepackage{graphicx}      
\usepackage{pgfplots}
\usepackage{natbib}        
\usepackage{amsmath,amssymb,mathtools}
\usepackage{algpseudocode}
\usepackage{empheq}
\usepackage{float} 
\DeclareMathOperator*{\argmin}{arg\,min}
\definecolor{myblue}{rgb}{.8, .8, 1}

\usepackage{pgfplots}
\pgfplotsset{compat=1.17} 
\usepackage{tikz}
\usetikzlibrary{plotmarks}
\usetikzlibrary{arrows.meta,bending,positioning,calc} 
\usetikzlibrary[shapes,arrows]
\usepgfplotslibrary{fillbetween}

\algnewcommand\algorithmicgiven{\textbf{Given:}}
\algnewcommand\Given{\item[\algorithmicgiven]}

\algnewcommand\algorithmicinitialise{\textbf{Initialise:}}
\algnewcommand\Initialise{\item[\algorithmicinitialise]}

\algnewcommand\algorithmicdefine{\textbf{Define:}}
\algnewcommand\Define{\item[\algorithmicdefine]}

\pgfplotsset{
    width=4.5cm, 
    height=3.5cm,
    axis line style={black},
    tick style={black},
    grid style={dashed, gray},
    label style={font=\tiny},
    tick label style={font=\tiny},
    legend style={font=\tiny},
    title style={font=\tiny, align=center, yshift=-0.2cm},
    legend cell align={left},
    legend style={
        fill opacity=0.8,
        draw opacity=1,
        text opacity=1,
        draw=gray,
        fill=white,
        anchor=south east, 
        at={(0.97,0.03)}, 
        column sep=0.1ex, 
        row sep=0.1ex, 
        inner sep=0.1ex, 
        legend image post style={line width=0.3mm, shorten <=1.5mm, shorten >=1.5mm} 
    },
}

\def\H{\mathcal{H}}
\def\R{\mathbb{R}}
\def\N{\mathbb{N}}
\def\prox{\operatorname{prox}}

\usepackage[acronym]{glossaries} 
\newacronym{PPA}{PPA}{Proximal Point Algorithm}
\newacronym{ccp}{ccp}{closed, convex and proper}
\newacronym{PVM}{PVM}{Proximal Value Method}
\newacronym{MPC}{MPC}{Model Predictive Control}
\begin{document}
\begin{frontmatter}

\title{A New Duality-Free Framework for Convex Optimisation with Superlinear Convergence and Effective Warm-Starting} 

\thanks[footnoteinfo]{This work was sponsored by the EPSRC.}

\author[First]{Michael Cummins} 
\author[First,Second]{Eric Kerrigan}  

\address[First]{Department of Electrical \& Electronic Engineering, Imperial College London, SW7 2AZ London, UK (email: m.cummins24@imperial.ac.uk)}
\address[Second]{Department of Aeronautics, Imperial College London, SW7 2AZ London, UK (email: e.kerrigan@imperial.ac.uk)}

\begin{abstract}                
Modern second-order solvers for convex optimisation, such as interior point methods, rely on primal-dual information and are difficult to warm-start, limiting their applicability in real-time control. We propose the \gls{PVM}, a duality‑free framework that reformulates the constrained problem as the unconstrained minimisation of a value function. The resulting problem always has a solution, yields a certificate of infeasibility and is amenable to warm‑starting. We develop a second‑order algorithm for Quadratic Programming based on the \gls{PPA} and semismooth Newton methods, and establish sufficient conditions for superlinear convergence to an arbitrarily small neighbourhood of the solution. Numerical experiments on a \gls{MPC} problem demonstrate competitive performance with state‑of‑the‑art solvers from a cold start and up to 70\% reduction in Newton iterations when warm starting.
\end{abstract}


\end{frontmatter}

\section{Introduction} 
Convex optimisation solvers are usually split into three categories: Interior--point methods \citep{goulart2024clarabel, coey2022solving, domahidi2013ecos, vandenberghe2010cvxopt}, operator splitting methods \citep{stellato2020osqp, garstka2021cosmo, o2021operator} and active set methods \citep{huangfu2018parallelizing, ferreau2014qpoases, arnstrom2022dual}. Interior point methods are generally regarded as the most robust option in many applications but are notoriously difficult to warm-start due to the central path \citep{wright1997primal}. Operator splitting and active set methods are effective when warm-started but also require a good guess for the dual variable, which can be difficult to infer. The goal of this paper is to make the first step towards a fourth option that can be effective from a warm start and remain competitive with interior point methods from a cold start.

Although dual variables can be useful in certain problems \citep{boyd2004convex}, they are not always beneficial when a sequence of similar convex problems need to be solved in quick succession such as real-time control \citep{gros2020linear}, trust-region subproblems \citep{nocedal2006numerical} or sequential convex programming \citep{malyuta2022convex}. In this sense, \gls{PVM} trades off the information gained from dual variables for better warm starting.

To achieve our desired trade-offs, we propose the Proximal Value Method (\gls{PVM}). Our method reformulates the problem as the unconstrained minimisation of a convex value function and solves it with the Proximal Point Algorithm (\gls{PPA}) \citep{rockafellar1976monotone}, using a semismooth Newton method for the arising subproblems \citep{qi1993nonsmooth}. We establish sufficient conditions for superlinear convergence and show that the resulting method matches the robustness of interior–point solvers while benefiting strongly from warm starts.

\textit{Related Work:}
Our approach is most closely related to FBStab \citep{liao2020fbstab}, a second-order method exhibiting superlinear convergence. Unlike \gls{PVM}, however, FBStab is restricted to QPs and relies on dual variable estimates for effective warm-starting. While \cite{bemporad2015quadratic} proposes a primal-dual active set method to recover dual warm starts from primal variables, it is similarly limited to QPs and lacks unboundedness detection, prohibiting its use as a general-purpose solver. In our case, \gls{PVM} generalises beyond QPs. The QP solver developed in our work is simply an application of \gls{PVM}.

\textit{Contributions:} Our contributions are as follows:
\begin{enumerate}
    \item We reformulate a broad class of convex optimisation problems into the unconstrained minimisation of a convex value function.
    \item We derive a tractable formulation of applying \gls{PPA} to the implicitly defined value function and provide sufficient conditions for superlinear convergence.
    \item We apply \gls{PVM} to a QP and derive semismooth Newton updates for solving the \gls{PVM} subproblems.
    \item We demonstrate \gls{PVM} significantly outperforms interior point solvers on a warm-started \gls{MPC} problem, achieving up to a 70\% reduction in Newton iterations.
\end{enumerate}

\textit{Notation:} We denote the non-negative orthant of $\R^n$ as $\R_+^n$. The Euclidean norm of a vector $x \in \R^n$ is denoted by $|x|$. The projection onto a closed convex set $\mathcal{C}$ is denoted by $\Pi_{\mathcal{C}}$. The subdifferential of a \gls{ccp} function $\pi:\R^n \to \R$ is denoted by $\partial \pi$. The sets of positive definite and positive semi-definite matrices in $\R^{n \times n}$ are denoted by $\mathbb{S}_{++}^n$ and $\mathbb{S}_+^n$, respectively. The distance from a vector $x \in \R^n$ to a closed convex set $\mathcal{C} \subseteq \R^n$ is denoted by $\operatorname{dist}(x,\mathcal{C}) := \inf_{y \in \mathcal{C}} |x - y|.$ We denote the componentwise operation of $\max\{\cdot,0\}$ to a vector $x\in\R^n$ as $(x)_+$.

The paper is structured as follows: Section~\ref{sec:prelim} presents the necessary preliminaries for \gls{PVM}. Section~\ref{sec:pvm} derives \gls{PVM} and provides sufficient conditions for superlinear convergence. In Section~\ref{sec:qp}, we specialise \gls{PVM} to a QP and show that these conditions can be verified. Section~\ref{sec:numerical} demonstrates the effectiveness of \gls{PVM} on a \gls{MPC} problem and we conclude the paper in section~\ref{sec:conclusion}.

\section{Preliminaries}\label{sec:prelim}
We review the semismooth Newton method as a means for solving the convex subproblems arising in \gls{PVM}, and \gls{PPA}, which serves as the core component of \gls{PVM}.

\subsection{Semismooth Newton}\label{sec:newton}
Due to the standard Newton method being restricted to functions that are twice differentiable, the semismooth Newton method was developed for a wider range of functions that satisfy a weaker property -- \emph{semismooth gradients} \cite[Chapter 1]{izmailov2014Newton}. Define $G:\R^n \to \R^m$ such that it is locally Lipschitz on a subset of its domain. The Clarke generalised Jacobian \citep{clarke1990optimization} may be defined as
\begin{align*}
    \partial_C G&(x) := \operatorname{co}\{ J \in \R^{m \times n} \: \mid \: \\&\exists\{x_k\} \subset D_G \: :
    \{x_k\} \to x, \: \{\nabla G(x_k)\} \to J\}, \nonumber
\end{align*}
where $D_G$ is a dense set of points where $G$ is differentiable. 
Moreover, for a differentiable function $g: \R^n \to \R$ whose gradient is (strongly) semismooth, its Clarke generalised Hessian is simply defined as 
\begin{align}\label{eq:gen_hess_def}
    \partial_C^2 g(x) := \partial_C (\nabla g)(x).
\end{align}
Now, we can define an iteration of the semismooth Newton method for minimising $g$ as
$$x_{j+1} = x_j + \rho_j\Delta x,$$
where the step direction $\Delta x$ satisfies 
$\H \Delta x = -\nabla g(x_j)$, $\quad \H \in \partial_C^2 g(x_j)$.
Similar to a pure Newton method, semismooth Newton is asymptotically quadratically convergent provided that $\partial_C^2g(x_j) \subset \mathcal{S}_{++}^n$ and $\{\rho_j\}$ is selected via a globalization mechanism. This powerful property will enable robust and accurate solutions to the subproblems arising in the proposed \gls{PVM} framework.

\subsection{The Proximal Point Algorithm}
Consider the optimisation problem $\min_x  h(x)$ and its corresponding optimality condition
\begin{equation}\label{eq:g_min}
0 \in \partial h(x),    
\end{equation}
where $h:\R^n \to \R$ is \gls{ccp}. \gls{PPA} computes a solution to~\eqref{eq:g_min} by solving a sequence of subproblems defined by the proximal operator
\begin{equation}
    x_{k+1} = \prox_{\sigma_k h}(x_k) := \argmin_x \: h(x) + \frac{1}{2\sigma_k}|x-x_k|^2. \notag
\end{equation}
A noteable property of \gls{PPA} is the strong convexity of each subproblem encountered at each iteration $k$. Moreover, each subproblem can be solved inexactly provided that the inexact iterate satisfies a summable error condition:
\begin{align}
\label{eq:A}
    |x_{k+1} - &\operatorname{prox}_{\sigma_k h}(x_k)|
\;\le\;
\varepsilon_k \min\{1, |x_{k+1} - x_k|^\varsigma\}\nonumber \\
&\sum_{k=0}^\infty \varepsilon_k < \infty, \quad \varsigma \geq 0. 
\tag{A}
\end{align}
For convergence to a solution of \eqref{eq:g_min}, it is sufficient to satisfy \eqref{eq:A} with $\varsigma = 0$. However, to discuss convergence rates, it is generally required to have $\varsigma \geq 1$ and that $\partial h$ satisfies a \textit{calmness} condition 
\begin{crit}[Calmness, \cite{rockafellar2021advances}]\label{as:error_bound}
    Suppose that \eqref{eq:g_min} has a solution set $S := \partial h^{-1}(0) \neq \emptyset$. Then, $\exists \: \vartheta > 0,\, \varphi > 0$ such that for some $y \in \R^n$ and $\bar y := \Pi_S(y)$,
    $$w \in \partial h(y), \:\: |w| < \vartheta \implies |y - \bar y| \,\leq\, \varphi|w|.$$
\end{crit}
Finally, since the proximal iterations in \gls{PVM} are implicitly defined, \eqref{eq:A} will allow us to construct a practical algorithm where each iteration is defined explicitly.

\section{The Proximal Value Method}\label{sec:pvm}
We first define the value function and study its basic structure. Then, we describe the idealised proximal point iteration on this scalar problem and extend it to a practical inexact scheme. The main results are Proposition~\ref{prop:r_star}, Lemma~\ref{lem:exact_outer}, Proposition~\ref{prop:two-step}, Theorem~\ref{thm:exact_outer} and the derivation of \gls{PVM}. All proofs are in Appendix~\ref{app:B}, with supporting results in Appendix~\ref{app:A}.

\subsection{The Value Function}

Consider the problem
\begin{equation}
\label{eq:QP}
\begin{aligned}
\min_{x,s} \quad & f(x), \\
\text{s.t.} \quad & A x  + s = b , \\
&s \in \mathcal{C},
\end{aligned}
\end{equation}
where $f: \R^n \to \R$ is \gls{ccp}, $\mathcal{C} \subseteq \R^m$ is a closed convex set, $x \in \mathbb{R}^n$, $s \in \R^m$, $A \in \mathbb{R}^{m \times n}$ and $b \in \mathbb{R}^m$.
We do not assume feasibility of \eqref{eq:QP}. An equivalent formulation of \eqref{eq:QP} can be obtained from the epigraphical reformulation
\begin{align*}
    \min_{x,s,t} \quad & t, \\
    \text{s.t.} \quad & f(x) \leq t,\\
    & Ax + s = b, \\
    & s \in \mathcal{C}.
\end{align*}
The goal of \gls{PVM} is to find the smallest $t \in\R$ such that the epigraphical constraints are feasible, which is equivalent to solving \eqref{eq:QP}. To this end, we define the \emph{residual vector field} $R:\mathbb{R}^n \times \mathbb{R}^m \times \mathbb{R} \to \mathbb{R}^{2+m}$ 
\begin{equation}
\label{eq:R-def}
R(x,s,t)
:=
\begin{pmatrix}
\max\{ f(x) - t, \; 0 \} \\
A x + s - b \\
\operatorname{dist}(s,\mathcal{C})
\end{pmatrix}.
\end{equation}

The associated merit function $r:\mathbb{R}^n \times \mathbb{R}^m \times \mathbb{R} \to \mathbb{R}_+$ is then defined by
\begin{align}
\label{eq:merit}
&r(x,s,t) \;\ :=\; \tfrac{1}{2} \, | R(x,s,t) |^2, \\
&= \frac{1}{2} \max\{ f(x) - t, \; 0 \}^2 + \underbrace{\frac{1}{2} (|Ax + s-b|^2 + \operatorname{dist}(s,\mathcal{C})^2)}_{r_0(x,s)}. \nonumber
\end{align}

We are now in a position to define the convex value function $r^\star:\mathbb{R} \to \mathbb{R}_+$ 
\begin{align}    
\label{eq:r_star_QP}
r^\star(t) \;\ :=\; \min_{x,s} r(x,s,t).
\end{align}
An equivalent problem to \eqref{eq:QP} can then be stated in terms of the value function as
\begin{equation}\label{eq:new_r_star_QP}
t^\star := \min_{t \in \R} \:t\quad \text{s.t.} \:\:  0 \in \partial r^\star(t).
\end{equation}

\begin{figure}
    \centering
    \input{plots/r_star}
    \caption{Value function $r^\star$ for two problems with objective $ f(x) = \frac{1}{2}x^2$ and equality constraint $x=1$. The inequality constraints are defined as $x \leq 0$ and $x \leq 2$ for the infeasible and feasible problems, respectively.}
    \label{fig:r_star}
\end{figure}
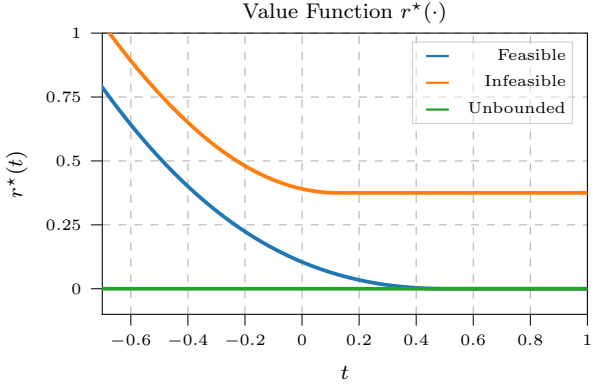

Figure~\ref{fig:r_star} provides an illustration for \eqref{eq:r_star_QP}. The interpretation is that if \eqref{eq:QP} is \emph{feasible}, then $r^\star(t)=0 \implies \exists \, (x,s) \: \text{ such that } A x + s = b,\; s \in \mathcal{C},\; f(x) \le t$ i.e., there exists a feasible point of cost at most $t$. If \eqref{eq:QP} is \emph{infeasible}, then $r^\star(t)>0$ for all $t$ and $r^\star(t)$ decreases to a positive constant as $t\to\infty$. However, both have a point where the slope flattens, and as a consequence, \eqref{eq:new_r_star_QP} always has a solution. The following formalises our discussion. 

\begin{prop}\label{prop:r_star}
    The following statements are true.
    \begin{enumerate}
        \item $r^\star$ is \gls{ccp} for all \gls{ccp} $f$ and closed convex $\mathcal{C}$.
        \item If \eqref{eq:QP} is feasible, $r^\star$ is strictly decreasing for all $t \leq t^\star$ and $r^\star(t) = 0$ for all $t \geq t^\star$.
        \item If \eqref{eq:QP} is infeasible, then $M := \min_{x,s} r_0(x,s) > 0$ and  $\lim_{t\to\infty} r^\star(t) = M$. Moreover, $\exists \tau \in \mathbb{R}$ such that $\partial r^\star(t) = 0$ and $r^\star(t) = M$ for all $t \geq \tau$ and $r^\star$ is strictly decreasing for all $t \leq \tau$. 
        \item If \eqref{eq:QP} is unbounded, then $r^\star(t) = 0$ for all $t \in \R$.
    \end{enumerate}
\end{prop}
\begin{rem}\label{rem:unbound}
    Although we don't detail a termination criterion for unboundedness, the theoretical basis is provided by Proposition~\ref{prop:r_star}. Future work could investigate methods for detecting such criteria.
\end{rem}
The remainder of this section presents iterative methods for solving \eqref{eq:new_r_star_QP} and provides sufficient conditions for superlinear convergence.

\subsection{Exact Proximal Iterations on the Value Function}

Although we do not solve \eqref{eq:new_r_star_QP} directly, we show that performing \gls{PPA} on \eqref{eq:r_star_QP} solves \eqref{eq:new_r_star_QP} under careful initialisation. Starting from an initial point $\hat t_0 \leq t^\star$, the \gls{PPA} iterations at time $k \in \N$ for solving \eqref{eq:new_r_star_QP} are
\begin{empheq}[box=\fbox]{align}
    \label{eq:prox_r_star}
    \hat t_{k+1} &= \prox_{\sigma_k r^\star}(\hat t_k). 
\end{empheq}

The following formalises the fact that \eqref{eq:prox_r_star} solves~\eqref{eq:new_r_star_QP} provided an appropriate initialisation.
\begin{lem}[Exact Convergence]\label{lem:exact_outer}
    Let Proposition~\ref{prop:r_star} hold and let the sequence $\{\hat t_k\}$ be generated by \eqref{eq:prox_r_star} with $\hat t_0 \leq t^\star$. Then,
    \begin{enumerate}
        \item $\{\hat t_k\}$ is monotonically non-decreasing and $\hat t_k \leq t^\star$ for all $k \in \N$
        \item $\{\hat t_k\}$ converges to $t^\star$. In particular, $\hat t_k \uparrow t^\star$
    \end{enumerate} 
\end{lem}

When discussing the convergence of \gls{PVM} throughout the paper, we will typically assume that $\hat{t}_0 \leq t^\star$ \footnote{This is generally possible in practice given problem knowledge such as norm--based constraints or the lower boundedness of $f$.}. In the case where $\hat{t}_0 > t^\star$, an adaptive lower-bound search can be used to compute a lower bound for $t^\star$, from which we can converge with superlinear convergence. Developing a mature solution to this problem could be the subject of future work.

Although \eqref{eq:prox_r_star} provides useful insight, it does not directly lead to a practical algorithm, since $r^\star$ is only defined implicitly. The rest of this section develops an iterative procedure that avoids this difficulty.

\subsection{Inexact Proximal Iterations on the Value Function}
First, we make the following assumption.
\begin{assum}\label{as:attainment_minimum}
    The merit function \eqref{eq:merit} attains its minimum in $(x,s)$ for a fixed $t \in \R$.
\end{assum}

A common sufficient condition for Assumption~\ref{as:attainment_minimum} is that~\eqref{eq:merit} is level-bounded. However, we later show that when applying \gls{PVM} to a QP or LP, Assumption~\ref{as:attainment_minimum} is satisfied immediately regardless of level-boundedness.
From this assumption, the following iterations are equivalent to \eqref{eq:prox_r_star} with respect to $\hat t$:
\begin{align}
\label{eq:Gk-def}
(x_{k+1}, s_{k+1}, &\hat t_{k+1})
\in 
\arg\min_{x,s,t} \: r(x,s,t) + \frac{1}{2\sigma_k} (t - \hat t_k)^2. 
\end{align}
This is the ``partial'' version of \gls{PPA} presented in \citet{bertsekas1994partial, ha1990generalization}. That is, we \emph{regularise only~$t$}, leaving $(x,s)$ unregularised. To this end, \cite{bertsekas1994partial} showed that \eqref{eq:Gk-def} is equivalent to a two-step update scheme. The scheme at iteration $k$ is defined as
\begin{empheq}[box=\fbox]{align}
    (x_{k+\frac{1}{2}}, s_{k+\frac{1}{2}}) &\in \argmin_{x,s} r(x,s,\hat t_k), \label{eq:Pa} \\
    (x_{k+1}, s_{k+1}, \hat t_{k+1}) &= \prox_{\sigma_k r}(x_{k+\frac{1}{2}}, s_{k+\frac{1}{2}}, \hat t_k). \label{eq:Pb} 
\end{empheq}

Therefore, solving \eqref{eq:Pa} and then \eqref{eq:Pb} is identical to solving~\eqref{eq:prox_r_star} with respect to the variable $\hat t$. In practice, we may not solve any of these subproblems exactly. This is well defined in terms of the \gls{PPA} method, where we initialise a sequence with $t_0 = \hat t_0$ and subsequent iterates are computed approximately with respect to the summable-error condition \eqref{eq:A}. Consider the stage objective of the proximal operator in~\eqref{eq:Pb} 
\begin{align}
    \phi_k(x,s,t)&:=r(x,s,t) \:+ \label{eq:Pb_stage}\\
    &\tfrac{1}{2\sigma_k}(|x-x_{k+\frac{1}{2}}|^2+|s-s_{k+\frac{1}{2}}|^2+(t-t_k)^2). \nonumber
\end{align}

The following details how to construct a sequence $\{t_k\}$ that satisfies \eqref{eq:A}, from only performing inexact evaluations of~\eqref{eq:Pa} and \eqref{eq:Pb}. 

\begin{prop}\label{prop:two-step}
Fix $k$ and let the iterates $(x_{k + \frac{1}{2}}, s_{k + \frac{1}{2}})$ and $(x_{k+1},s_{k+1}, t_{k+1})$ be bounded. Additionally, let $\mathcal{S}_k := \argmin_{x,s} r(x,s,t_k)$ and define  
\begin{align*}
\delta_k &:= \nabla \phi_k(x_{k+1},s_{k+1}, t_{k+1}).
\end{align*}
Then, the following bound holds for all $k \in \N$:
$$
|t_{k+1}-\operatorname{prox}_{\sigma_k r^\star}(t_k)| \;\le\; \sigma_k
| \delta_k | \,+\, \operatorname{dist}((x_{k+\frac{1}{2}},s_{k+\frac{1}{2}}),\mathcal{S}_k).
$$
\end{prop}

    Deriving an explicit upper bound for the distance term can be achieved through the use of \emph{error bounds} or \emph{growth conditions} \citep{pang1997error}. Such bounds are typically represented in terms of the gradients of the problem-defining objective. 
We now have some of the necessary ingredients to construct practical algorithms using \gls{PVM}. 

\begin{thm}[Inexact Convergence]\label{thm:exact_outer}
    Let $\{\hat t_k\}$ be generated by \eqref{eq:prox_r_star}, and let $\{t_k\}$ be any sequence satisfying \eqref{eq:A} with $t_0=\hat t_0$ and $\varsigma \geq 0$. Suppose Lemma~\ref{lem:exact_outer} holds.
    Then, $\{t_k\}$ is well defined and converges to a limit $\tau$ with
    \[
    \tau \in [t^\star, t^\star + \lambda],
    \]
    where $\lambda := \sum_{k=1}^\infty \varepsilon_k\min\{1,|t_{k+1} - t_k|^\varsigma\}$ can be made arbitrarily small. Additionally, if $\partial r^\star$ satisfies Criterion~\ref{as:error_bound}, $\sigma_k \uparrow \infty$ and $\varsigma \geq 1$, then the convergence rate is superlinear. If $\sigma_k$ remains constant, the rate is linear.
\end{thm}

Although the overshoot $\lambda$ may appear to be inconvenient, we later show in our experiments that even when $\lambda$ is orders of magnitude larger than $t^\star$, \gls{PVM} still converges to $t^\star$ within a tolerance of $10^{-8}$. We hypothesise that the overshoot bound $\lambda$ is highly conservative, and proving a tighter bound could be the subject of future work. 

\section{Solving a Quadratic Program}\label{sec:qp}
In this section, we specialise \gls{PVM} to the case where \eqref{eq:QP} is either a QP or an LP. Specifically, we have
\begin{align}\label{eq:QP_specification}
    \min_{x,s} \quad & q(x) := \tfrac{1}{2} x^\top Q x + p^\top x, \nonumber\\
    \text{s.t.} \quad & A x  + s = b , \\
    &s \geq 0, \nonumber
\end{align}
where $Q \in \mathbb{S}_+^n$. Since $\mathcal{C} = \R^M_+$ we can then write 
$$\operatorname{dist}(s,\mathcal{C}) \,=\, |(I - \Pi_{\R^m_+})s| \,=\, |(-s)_+| .$$
To this end, \eqref{eq:merit} is differentiable and has strongly semismooth gradients. This section describes a procedure for applying \gls{PVM} to a specific convex optimisation problem, such as \eqref{eq:QP_specification}. First, we provide an explicit way to control the overshoot in Theorem~\ref{thm:exact_outer}, using Proposition~\ref{prop:two-step}. Next, we describe how to apply the semismooth Newton method to subproblems \eqref{eq:Pa} and \eqref{eq:Pb}, with specific attention given to the linear systems. Finally, we prove that Assumption~\ref{as:attainment_minimum} and Criterion~\ref{as:error_bound} hold for the specifications $\mathcal{C} = \R_+$ and $f(x) = q(x)$, enabling us to prove superlinear convergence. The main results are Proposition~\ref{prop:inner},  Theorem~\ref{thm:QP_convergence} and the application of \gls{PVM} to  \eqref{eq:QP_specification}. All proofs are in Appendix~\ref{app:B}. 

\subsection{Controlling the Overshoot}
First, we make the following assumption.

\begin{assum}[Quadratic Growth]\label{as:qg}
    Let $(x_{k+\frac{1}{2}}, s_{k+\frac{1}{2}})$ be bounded for any given $k \in \N$. Also, $\exists \, \kappa > 0$ such that 
    $$\operatorname{dist}((x_{k+\frac{1}{2}}, s_{k+\frac{1}{2}}), \mathcal{S}_k) \leq \kappa |\nabla r(x_{k+\frac{1}{2}}, s_{k+\frac{1}{2}}, t_k)|,$$
    holds for all $t_k \in \R$.
\end{assum}
    Quadratic growth is a weaker condition than strong convexity and holds for many convex functions with non-unique minima \citep{li1995error, pang1997error}. 
    For example, if $Q = 0$, then Assumption~\ref{as:qg} automatically holds since \eqref{eq:merit} is convex piecewise-linear-quadratic \citep{li1995error}.
Now, we use Proposition~\ref{prop:two-step} to drive the overshoot in Theorem~\ref{thm:exact_outer} to~0. Let $\delta_k^\alpha \in \R_+$ and $\delta_k^\beta \in \R_+$ be the gradient tolerances used at iteration $k \in \N$ to terminate the iterative methods (possibly semismooth Newton), minimising \eqref{eq:Pa} and \eqref{eq:Pb}, respectively. Note that by setting $\sigma_k = 1/\sqrt{\delta_k^\beta}$ in \eqref{eq:Pb}, we have from Assumption~\ref{as:qg} and Proposition~\ref{prop:two-step} that
$$|t_{k+1} - \prox_{\sigma_k r^\star}(t_{k})| \leq \sqrt{\delta_k^\beta} + \kappa \delta_k^\alpha.$$

Therefore, we have full control over the overshoot and can enforce it to be arbitrarily small. Moreover, we can guarantee \eqref{eq:A} for any $\varsigma \geq 0$ provided that we choose appropriate sequences $\{\delta_k^\beta\}$ and $\{\delta_k^\alpha\}$.

\subsection{Semi-Smooth Newton}
Suppose for a moment that we are only concerned with solving \eqref{eq:Pa}. First, we define the utility functions
\begin{align}
    r_q(x,t) := &\max\{q(x) - t, 0\} \label{eq:rq},\\
    \xi_q(x,t) := &\begin{cases}
        1, \:&q(x) > t, \\
        0, \:&q(x) < t, \\
        \theta \in  [0,1], &q(x) = t,
    \end{cases} \label{eq:xi_q} \\
    \xi_s(s_i) := &\begin{cases}
        1, \:&s_i < 0, \\
        0, \:&s_i > 0, \\
        \theta \in  [0,1], &s_i = 0.
    \end{cases} \label{eq:xi_s}
\end{align}
Then, defining $\Gamma(s) := \operatorname{diag}(\xi_s(s_1),\dots,\xi_s(s_m)) $ and $v := Qx + p$, the generalised Hessians of \eqref{eq:merit} with respect to $x$ and $s$ are 
\begin{align}
    \partial_{C,x}^2r(x,s,t) = \:&A^\top A + Qr_q(x,t_k) + \xi_q(x,t) vv^\top, \label{eq:x_hess_r}\\
    \partial_{C,s}^2r(x,s,t) = \:&I + \Gamma(s), \label{eq:s_hess_r}
\end{align}
with the full generalised Hessian of \eqref{eq:merit} being  
\begin{align}
    \partial_C^2 r(x, s, t) = \left\{
    \begin{pmatrix}
        H_x & A^\top \\ A & H_s 
    \end{pmatrix}  : 
     H_x \in \eqref{eq:x_hess_r}, \: H_s \in \eqref{eq:s_hess_r}
    \right\}. 
\end{align}
A semismooth Newton iteration will require an arbitrary matrix from the generalised Hessian and then solve the resulting linear system to compute a step direction. 
To this end, we define $\bar v := (v,0)$, restrict $\theta = 1$ in \eqref{eq:xi_q} and~\eqref{eq:xi_s}, and consider an element of the generalised Hessian at Newton iteration $j \in \N$ in the form 
\begin{equation}
    \H(x_j,s_j,t_k) + \xi_q(x_j,t_k)\bar v \bar v^\top \in \partial_C^2 r(x_j,s_j,t_k).  
\end{equation}
This reformulation will allow us to solve the linear systems using a sparse factorisation of $\H(x_j,s_j,t_k)$ and a Sherman-Morrison update \citep{petersen2008matrix} to handle the dense rank-$1$ term $\xi_q(x_j,t_k)\bar v \bar v^\top$.
If we instead want to minimise~\eqref{eq:Pb_stage}, define $\bar w = (\bar v, -1)$ and observe an element of the generalised Hessian 
\begin{align}
     \begin{pmatrix}
        \H(x_j,s_j,t) + \sigma_k^{-1}I & 0\\ 0 & \sigma_k^{-1}
    \end{pmatrix} + \xi_q(x_j,&t)  \bar w \bar w^\top \\ &\in \: \partial_C^2\phi_k(x_j,s_j,t). \nonumber
\end{align}
Therefore, whether we are solving subproblem \eqref{eq:Pa} or~\eqref{eq:Pb}, we factorise the same matrix $\H(x_j,s_j,t_j)$ at each Newton step. Moreover, $\H(x_j,s_j,t_j)$ can only have two distinct sparsity patterns depending on whether $r_q(x_j,t_j) \geq 0$. This observation supports our use of allocation-free sparse matrix factorisation tools, leading to significantly lower factorisation complexity and enabling real-time execution.

One remaining issue is the possible singularity of the linear systems arising in subproblem \eqref{eq:Pa}, since the merit function $r$ may not be strongly convex. To this end, we consider the augmented Hessian at iteration $j \in \N$
$$\bar \H(x_j,s_j,t_k) := \H(x_j,s_j,t_k) + \mu_j I,$$
where $\mu_j := |\nabla r(x_j,s_j,t_k)|$. Provided Assumption~\ref{as:qg} holds and $\H(x_j,s_j,t_k) \succeq 0$, asymptotic superlinear convergence to the set of solutions $\mathcal{S}_k$ is guaranteed \cite[Chapter 3]{ulbrich2011semismooth}.

\subsection{Convergence Analysis}
First, we study the convergence of the subproblems \eqref{eq:Pa} and \eqref{eq:Pb} under semi-smooth Newton iterations. 

\begin{prop}[Inner Convergence]\label{prop:inner}
    Consider \eqref{eq:Pa} and \eqref{eq:Pb} at a fixed outer iteration $k \in \N$, with solutions $(x^\star_{k+\frac{1}{2}}, s^\star_{k+\frac{1}{2}})^\alpha$ and $(x^\star_{k+1}, s^\star_{k+1}, t^\star_{k+1})^\beta$, respectively. Then, the sequences generated by semi-smooth Newton with globalisation \cite[Algorithm 3.1]{nocedal2006numerical} are well-defined and satisfy 
    \begin{enumerate}
        \item $\{(x_j, s_j)^\alpha\}_j \to (x^\star_{k+\frac{1}{2}}, s^\star_{k+\frac{1}{2}})^\alpha$,
        \item $\{(x_j, s_j, t_j)^\beta\}_j \to (x^\star_{k+1}, s^\star_{k+1}, t^\star_{k+1})^\beta$,
    \end{enumerate} 
    respectively. Moreover, the asymptotic rate of convergence is superlinear for (1) and quadratic for (2). 
\end{prop}

The exceptional convergence rate and scale invariance of the semismooth Newton method will allow us to set high gradient tolerances during inner subproblems and drive the overshoot in Theorem~\ref{thm:exact_outer} to an arbitrarily small value. Finally, we verify the sufficient conditions required for the convergence of \gls{PVM} when $\mathcal{C} = \R_+^m$ and $f(x) = q(x)$.

\begin{thm}[QP Convergence]\label{thm:QP_convergence}
    Let \eqref{eq:QP} satisfy \eqref{eq:QP_specification}. Then, the following statements are true:
    \begin{enumerate}
        \item The merit function \eqref{eq:merit} attains its minimum in $(x,s)$ for any $t\in\R$.
        \item The operator $\partial r^\star$ satisfies Criterion~\ref{as:error_bound}.
        \item \gls{PVM} initialised with $t_0 \leq t^\star$ will converge to a neighbourhood of $t^\star$ at a superlinear rate if $\sigma_k \uparrow \infty$, or at a linear rate if $\sigma_k$ is constant.
    \end{enumerate}
\end{thm}

Although \gls{PVM} requires certain assumptions to verify convergence, we argue that they can be easy to verify in practice and \gls{PVM} may therefore emerge as a suitable alternative to interior point methods in certain problems. 

\section{Numerical Evaluation}\label{sec:numerical}
Our primary evaluation for PVM will be a linear-quadratic Model Predictive Control (\gls{MPC}) problem. The \gls{MPC} problem with state and input variables at time $i$, denoted as $(z_i, u_i)$, is defined as 
\begin{align}\label{eq:mpc}
    \min_{z, u}&\:\: \sum_{i=1}^H z_{i}^\top Q z_i + u_{i}^\top R u_i \\
    \text{s.t.} &\:\: z_0 = \hat z_0, \:\: z_{i+1} = A_{\text{sys}}z_i + B_{\text{sys}}u_i, \quad i = 1, \dots, H, \nonumber \\
    &\:\: |u_i|_{\infty} \leq u_{\text{max}}, \quad i=1,\dots,H \nonumber \\
    &\:\: |z_i|_{\infty} \leq z_{\text{max}}, \quad i=1,\dots,H \nonumber
\end{align}

where $\hat z_0$ is the initial state estimate, $H=20$, $Q = 100I_{3\times3}$, $R = 0.01I_{3\times3}$, $u_{\text{max}} = 0.01$, $z_\text{max} = 1$ and the system matrices are 
$$A_{\text{sys}} = \begin{pmatrix}
    1.01 & 0.01 & 0 \\ 0.01 & 1.01 & 0.01 \\ 0 & 0.01 & 1.01 
\end{pmatrix}, \quad B_{\text{sys}} = I_{3 \times 3}.$$
We convert \eqref{eq:mpc} into a QP of the form \eqref{eq:QP_specification} by letting the primal variable of \eqref{eq:QP_specification} be $x = (z_1,\dots,z_H,u_1,\dots,u_H)$.
The constraint matrices in \eqref{eq:QP} are then tall and sparse with $\mathcal{C} = \R_+^m$ and the number of decision variables $(x,s)$ is 495. We developed a prototype QP solver in Julia \citep{bezanson2017julia} and compare against three state-of-the-art open-source solvers: Clarabel, ECOS and Hypatia.

The metric used for comparison is the total number of Newton steps as the primary computational bottleneck in both \gls{PVM} and interior point solvers is the matrix factorisation at each step. For \gls{PVM}, the number of Newton iterations is the cumulative number of Newton steps taken in all subproblems from the beginning to convergence. Additionally, when discussing accuracy of a solution, we refer to $|x_c - x^\star|_{\infty}$, where $x^\star$ is the \gls{PVM} solution and $x_c$ is the Clarabel solution.

\subsection{Parameter Selection and Termination}

We describe a scheduling for the parameters $\{\sigma_k\}$, $\{\delta_k^\alpha\}$, and $\{\delta_k^\beta\}$ that should produce a huge overshoot; however, in practice, it produces a sequence satisfying $t_k \uparrow t^\star$. Consider the scheduling
$$ \sigma_{k+1} = \max\{1/\sqrt{\delta_k^\beta}, \sigma_{k}\}, \quad \delta_{k+1}^\beta = \delta_{k+1}^\alpha = \delta_k^\beta / 10,$$
with initial parameters $\sigma_0 = 10^4$ and $\delta_0^\beta = 10^{-2}$. Due to Theorem~\ref{thm:QP_convergence}, we therefore terminate when both
$$|\nabla_t r(x_{k+\frac{1}{2}}, s_{k+\frac{1}{2}}, t_k)| \leq \epsilon_{\text{opt}}, \quad \delta_k^\beta \leq \epsilon_{\text{opt}},$$
with optimality tolerance $\epsilon_{\text{opt}} = 10^{-4}$. We deem the solution feasible if $r(x_{k+\frac{1}{2}}, s_{k+\frac{1}{2}}, t_k) \leq \epsilon_{\text{con}}$, with constraint tolerance $\epsilon_{\text{con}} = 10^{-8}$.

\subsection{Warm-Starting}
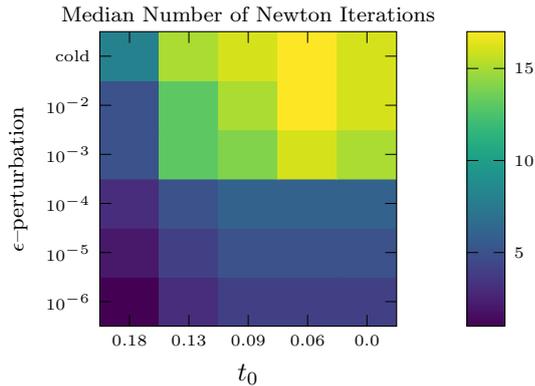
\begin{figure}[tb]
    \centering
    \input{plots/median_matrix}
    \caption{Median number of Newton iterations required to reach convergence when perturbing $t$ and $(x,s)$ warm start. All problems are solved to within $10^{-8}$ accuracy of the true solution.}
    \label{fig:heatmaps}
\end{figure}
We evaluate the efficacy of warm starting by varying the initial distances from the primal solution and from $t_0$ to $t^\star$. We scale the objective of $\eqref{eq:mpc}$ such that $t^\star = 0.1819$ and consider different $t_0$ in the range $[0,0.1819]$. We compute the optimal solution to \eqref{eq:mpc} using Clarabel, denoted as $(x^\star,s^\star)$ and then set $x_0 = x^\star + \bar \epsilon$, where $\bar \epsilon$ is a random vector with each element sampled from the uniform distribution defined on the range $[-\epsilon,\epsilon]$. Additionally, we show the number of iterations from a cold start where $x_0$ is simply the zero vector. For the slack variable, we set $s_0 = b - Ax_0$. Figure~\ref{fig:heatmaps} illustrates our results with the value of $\epsilon$ on the y-axis and $t_0$ on the x-axis. The iterations required for Clarabel, ECOS and Hypatia are $9$, $14$ and $16$, respectively

Note \gls{PVM} matches ECOS and Hypatia from a cold start, highlighting the robustness of our method and its superlinear convergence. Effective warm-starting yields a 70\% reduction in iterations compared to Clarabel, requiring only a single step in the near-solution regime. Solutions agree within an accuracy $10^{-8}$, eluding to our hypothesis that the theoretical overshoot in Theorem~\ref{thm:exact_outer} is highly conservative.

\subsection{Recovering Feasibility}
\begin{figure}[tb]
    \centering
    \input{plots/feasible}
    \caption{Statistics for total number of Newton iterations required to reach convergence when perturbing $t$ and $(x,s)$ warm start.}
    \label{fig:recover}
\end{figure}
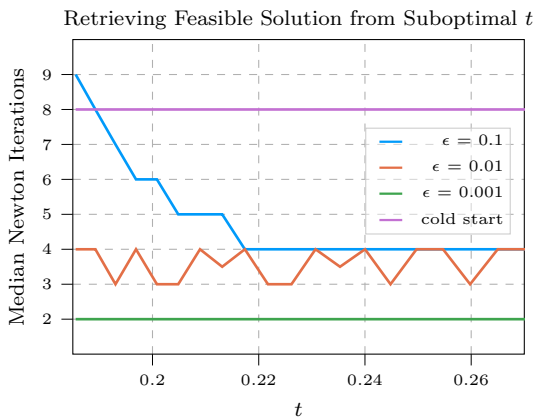
An important property of \gls{PVM} is the ability to recover feasible solutions from a suboptimal $t_0 > t^\star$. We generated 20 different increasing values for $t_0$ that are all larger than $t^\star$. We compute the solution $(x^\star(t_0), s^\star(t_0)) \in \argmin_{x,s} r(x,s,t_0)$ and similarly to the warm-starting evaluation, perturb the initial $x_0$ and report the number of iterations to reach convergence. Results are presented in Fig.~\ref{fig:recover}. As expected, \gls{PVM} performs very well at computing feasible solutions even from a cold start. Future work could explore the use of \gls{PVM} in real-time optimal control where feasible solutions are prioritised.

\subsection{Infeasiblity Detection}
\begin{table}[tb]
\caption{Statistics for number of iterations required to detect infeasibility of \eqref{eq:QP}.}
    \centering
    \begin{tabular*}{0.8\linewidth}{@{\extracolsep{\fill}}lccc@{}}
        \hline
        \textbf{Solver} & \textbf{Min } & \textbf{Max} & \textbf{Median } \\
        \hline
        ECOS     & 4.0 & 6.0  & 5.0  \\
        Hypatia  & 4.0 & 27.0 & 14.0 \\
        Clarabel & 5.0 & 7.0  & 6.0  \\
        PVM      & 5.0 & 9.0  & 7.0  \\
        \hline
    \end{tabular*}
    
    \label{tab:infeasible}
\end{table}
We construct infeasible problems from \eqref{eq:mpc} by appropriately selecting values of $\hat s_0$, $u_{\text{max}}$ and $s_{\text{max}}$. First, we restrict $\hat s_0 = (1,1,1)$ and then select the other parameters from a range of values $u_{\text{max}} \in \{-0.01, -0.02, \dots -0.09\}$ and $s_{\text{max}} \in \{0.1, 0.15, \dots, 0.5\}$, respectively. This leads to a total of 72 infeasible problems. For each solver, we report the minimum, maximum and median number of iterations required for convergence, respectively. Additionally, all \gls{PVM} solves are operated from a cold start. Results are displayed in Table~\ref{tab:infeasible}. \gls{PVM} remains competitive against ECOS and Clarabel and is significantly more efficient than Hypatia. Following this observation, future work could investigate the use of \gls{PVM} as an inner solver in mixed integer problems, where infeasible problems are encountered frequently.

\section{Conclusions and Future Work}\label{sec:conclusion}
We introduced \gls{PVM}, a new duality-free framework for convex optimization that converts the constrained problem into the unconstrained minimisation of a convex value function. We derived a tractable formulation for solving this problem and specialised our method to a QP. Sufficient conditions for superlinear convergence were presented and numerical evaluation demonstrated the ability to achieve highly accurate solutions in very few Newton iterations, particularly when warm-started. Although the proposed framework is still in its infancy, we believe there is significant promise in the current results. Moreover, due to the huge success of \gls{PPA} in nonconvex optimisation, we believe \gls{PVM} will also extend nicely to non-convex problems upon further investigation.

There are many directions for future work. Firstly, one could apply \gls{PVM} to derive algorithms for a wider range of convex problems and also to establish convergence guarantees in the non-convex domain. Second, the overshoot bound could be investigated further in an effort to derive an optimal bound. Finally, a more general strategy needs to be developed for handling scenarios where $t_0 > t^\star$ and detecting unboundedness (Remark~\ref{rem:unbound}).

\bibliography{refs}             

\appendix 
\section{Supporting Results}
\label{app:A}

\begin{lem}\label{lem:quadratic-growth-rstar} 

Let $r^\star(\cdot)$ be defined as in~\eqref{eq:r_star_QP} and \eqref{eq:QP_specification} and let $M = \min_{x,s} r_0(x,s)$. Then, $\exists \: a > 0$ such that
\[
r^\star(t) \;\ge\; a\,(t - t^\star)^2 + M,
\quad \forall\, t \leq t^\star.
\]
\end{lem}

\textit{Proof:} Let
$
X := \argmin_{x,s} r_0(x,s)
$
be the set of points closest to the (empty) feasible set of \eqref{eq:QP}. To this end, $X$ is non-empty \citep{belousov2002frank} and polyhedral \citep{mangasarian1988simple}. Therefore, $X$ is governed by a system of linear inequalities that always has a solution. Define the epigraphical feasible set
\[
S := \{(x,s,t) \in \mathbb{R}^{n+m+1} \mid (x,s) \in X,\ q(x) \le t\}.
\]
Note that $r(x,s,t) = 0$ if and only if $(x,s,t) \in S$. Let $(y,u,\tau) \in S$. Then $(y,u) \in X$ and $q(y) \le \tau$. By definition of $t^\star$ as the minimum of $q$ over $X$, we have $q(y) \ge t^\star$, hence $\tau \ge q(y) \ge t^\star$. Fix any $t \leq t^\star$ and let $(x^\star(t),s^\star(t),t) \in \argmin_{x,s} r(x,s,t)$. Then for every $(y,u,\tau) \in S$,
\[
|(x^\star(t),s^\star(t),t) - (y,u,\tau)|
\;\ge\; |t - \tau|
\;\ge\; t^\star - t.
\]

In particular, for the minimiser $(x^\star(t),s^\star(t))$ we have
\begin{equation}\label{eq:dist-lower-minimiser}
\operatorname{dist}\bigl( (x^\star(t),s^\star(t),t), S \bigr)
\;\ge\; t^\star - t
\quad \forall\, t \leq t^\star.
\end{equation}

Now, note that the set $S$ is governed by convex quadratic inequalities. Moreover, for any $(x,s,t) \in S$, we also have that $(x,s,t+\delta) \in S$ $\forall \, \delta \geq 0$; therefore, the system is non-singular. To this end, we can invoke \cite[Corollary 12]{pang1997error} and deduce $\exists \: \kappa > 0$ such that 
$$\operatorname{dist}((x,s,t),S) \leq \kappa \max\{\operatorname{dist}((x,s),X),\,(q(x) - t)_+\},$$

for all  $(x,s,t) \in \R^{n+m+1}$. From \cite[Theorem 2.7]{li1995error}, we also know that $\exists\: \gamma > 0$ such that 
$$\operatorname{dist}((x,s),X) \leq \gamma \sqrt{r_0(x,s) - M.}$$

To this end, $\exists\:c>0$ such that 
\begin{align}
 \operatorname{dist}((x,s,t),S) \leq &\:c \sqrt{(q(x) - t)_+^2 + r_0(x,s) - M}, \nonumber\\
 =&\:c\sqrt{r(x,s,t) - M},\nonumber
\end{align}
for all  $(x,s,t) \in \R^{n+m+1}$. Combining with \eqref{eq:dist-lower-minimiser} yields 
\[
r^\star(t)
\;\ge\; \frac{1}{c^2}\,(t^\star - t)^2 + M.
\]
for all $t \leq t^\star$. Setting $a := \frac{1}{c^2} > 0$ proves the claim. \qed

\section{Proofs of Main Results}\label{app:B}

\textit{Proof of Propostion~\ref{prop:r_star}:} We know that \eqref{eq:merit} is convex and therefore \eqref{eq:r_star_QP} is also convex. To show \eqref{eq:r_star_QP} is proper, let $\bar x \in \mathcal{C}$ and note that $r(\bar x, t) < \infty$ for any t. Then, $0 \leq r^\star(t) \leq r(\bar x, t) \: \forall \, t \in \R$ which implies that \eqref{eq:r_star_QP} is proper. To this end, since \eqref{eq:r_star_QP} is finite and convex, it is also continuous \cite[Theorem 10.1]{rockafellarconvex} and therefore closed \cite[Theorem 7.1]{rockafellarconvex}. 

2) follows directly from 1) and the fact that $r^\star$ is non-increasing over its domain \cite[Theorem 1]{sequoia}. 

To prove 3) we first note that $\lim_{t \to \infty} \max\{f(x) - t\} = 0$ which concludes $M = \lim_{t \to \infty} r^\star(t)$. Now, let $(\bar x, \bar s) \in \argmin r_0(x,s)$ and $\bar \tau  = q(\bar x)$. Then, $r(\bar x, \bar s, t) = r^\star(t) = M$ for all $t \geq \bar \tau$. Moreover, due to $r^\star$ being non-increasing and convex, there must exist a $\tau \leq \bar \tau$ such that $r^\star$ is strictly decreasing for $t \leq \tau$. 

4) follows from \cite{sequoia}. \qed

\textit{Proof of Lemma~\ref{lem:exact_outer}:} For convex $r^\star$, $t=\hat t_{k+1}$ solves \eqref{eq:prox_r_star} if and only if
\begin{equation}
\label{eq:FOC}
0 \in \partial r^\star(\hat t_{k+1})
+ \frac{1}{\sigma_k}(\hat t_{k+1} - \hat t_k),
\end{equation}
i.e.\ there exists $g_{k+1} \in \partial r^\star(\hat t_{k+1})$ such that
\begin{equation}
\label{eq:prox-eq}
0 = g_{k+1} + \frac{1}{\sigma_k}(\hat t_{k+1} - \hat t_k)
\:\Longleftrightarrow\:
\hat t_{k+1}
= \hat t_k - \sigma_k g_{k+1}.
\end{equation}

By Proposition~\ref{prop:r_star},
$r^\star$ is strictly decreasing on $(-\infty,t^\star]$ and becomes flat (constant) on $[t^\star,\infty)$.
Hence, if $u < t^\star$, then $\partial r^\star(u) < 0$, and if $u \geq t^\star$, then $\partial r^\star(u)=0$.

Inductively, $\hat t_k \le t^\star$ for all $k$.
If $\hat t_k \le t^\star$ and $\hat t_{k+1} > t^\star$, then $g_{k+1}=0$ by flatness, which by \eqref{eq:prox-eq} gives $\hat t_{k+1}=\hat t_k$, a contradiction. Thus $\hat t_{k+1}\le t^\star$.
Moreover, with $g_{k+1}\le 0$ we have $\hat t_{k+1}=\hat t_k-\sigma_k g_{k+1}\ge \hat t_k$, proving (a).

Since $\{\hat t_k\}$ is non-decreasing and bounded above, it converges to some $\bar t\le t^\star$.
From \eqref{eq:prox-eq} and \cite{rockafellar1976monotone},
$g_{k+1}=-(\hat t_{k+1}-\hat t_k)/\sigma_k\to 0$ as $k\to\infty$.
Closedness of $\operatorname{graph}\partial r^\star$ yields $0\in\partial r^\star(\bar t)$, i.e.\ $\bar t\in\arg\min r^\star=[t^\star,\infty)$.
Hence $\bar t=t^\star$, and $\hat t_k\uparrow t^\star$.
\qed

\textit{Proof of Propostion~\ref{prop:two-step}:} Define the shorthand $z := (x,s)$ such that $(\bar z_{k+1}, \bar t_{k+1}) := \prox_{\sigma_k r}(\bar z_{k+\frac{1}{2}}, t_k)$ and let $\bar z_{k+\frac{1}{2}} \in \mathcal{S}_k$. By the strong convexity of \eqref{eq:Pb_stage} and the corresponding PL inequality, we have  

    $$|\,(z_{k+1}, t_{k+1}) - \prox_{\sigma_k r}(z_{k + \frac{1}{2}},t_k)\,| \leq \sigma_k |\delta_k|.$$
    Moreover, we known that $\bar t_{k+1} = \operatorname{prox}_{\sigma_k r^\star}(t_k)$ and therefore we can say 
    \begin{align*}
        |\,t_{k+1}&-\prox_{\sigma_k r^\star}(t_k)\,| \leq \:|\,(z_{k+1}, t_{k+1}) - \prox_{\sigma_k r}(\bar z_{k+\frac{1}{2}}, t_k)\,|,    \\
            &\leq \:  |\,(z_{k+1}, t_{k+1}) - \prox_{\sigma_k r}(z_{k+\frac{1}{2}}, t_k)\,| \\
            &\qquad\qquad\: + |\,\prox_{\sigma_k r}(z_{k+\frac{1}{2}}, t_k) - \prox_{\sigma_k r}(\bar z_{k+\frac{1}{2}}, t_k) \,|, \\
        &\leq \:  \sigma_k |\delta_k| + |\,\prox_{\sigma_k r}(z_{k+\frac{1}{2}}, t_k) - \prox_{\sigma_k r}(\bar z_{k+\frac{1}{2}}, t_k) \,|,
    \end{align*}
    for all $\bar z_{k+\frac{1}{2}} \in \mathcal{S}_k$. Finally, by the firm non-expansiveness of $\prox_{\sigma_k r}$ we can conclude
    \begin{align*}
        |\,t_{k+1}-\operatorname{prox}_{\sigma_k r^\star}(t_k)\,| &\leq \: \sigma_k |\delta_k| + |\,(z_{k + \frac{1}{2}},t_k) - (\bar z_{k + \frac{1}{2}},t_k) \,| \\
        &= \sigma_k |\delta_k| + \operatorname{dist}(z_{k+\frac{1}{2}}, \mathcal{S}_k). 
    \end{align*} 

    This concludes the proof. \hfill \qed

\textit{Proof of Theorem~\ref{thm:exact_outer}:} Let $\hat t_{k+1}:=\operatorname{prox}_{\sigma_k r^\star}(\hat t_k)$ with $\hat t_0=t_0$. Define $S_k:=\sum_{i=0}^{k-1}\varepsilon_i\min\{1,|t_{i+1} - t_i|^\varsigma\}$ for $k\ge1$ and $S_0:=0$.
Using (\normalfont A) and the firm non-expansiveness of $\operatorname{prox}_{\sigma_k r^\star}$ (in particular, $1$-Lipschitz),
\begin{align}
|t_{k+1} - \hat t_{k+1}| \le 
|\,t_{k+1}-\operatorname{prox}_{\sigma_k r^\star}(t_k)\,|
\; &+\\ |\,\operatorname{prox}_{\sigma_k r^\star}(t_k)-\operatorname{prox}_{\sigma_k r^\star}(\hat t_k)\,| \; \le \;
\varepsilon_k&\min\{1,|t_{k+1} - t_k|^\varsigma\} \nonumber \\
&\quad+ |t_k-\hat t_k|.   \nonumber 
\end{align}

Since $S_{k+1} = \varepsilon_k\min\{1,|t_{k+1} - t_k|^\varsigma\} + S_k$, we can also conclude $|t_k-\hat t_k|\le S_k$ for all $k$ by induction. 
From \cite{rockafellar1976monotone}, we know that under condition (\normalfont{A}) ${t_k}$ is well defined and converges to a value $\tau := \limsup_{k \to \infty} t_k$ such that $\tau \in [t^\star, \infty)$. Following from Lemma~\ref{lem:exact_outer} and the previous inequality we have that 
\begin{align}
    \tau = \limsup_{k\to\infty} t_k
    &\le \lim_{k\to\infty}\hat t_k + \limsup_{k\to\infty}|t_k-\hat t_k| \nonumber\\
    &\le t^\star + \lim_{k\to\infty} S_k
    = t^\star+\lambda.    \nonumber
\end{align}

\noindent
Then, we can say $\tau \in [t^\star, \infty) \: \implies \: \tau \in [t^\star, t^\star + \lambda]$. 
Finally, since the operator $\partial r^\star$ satisfies Criterion~\ref{as:error_bound}, the result of linear and superlinear convergence follows directly from \cite[Theorem 2.2]{rockafellar2021advances}. \hfill \qed

\textit{Proof of Theorem~\ref{thm:QP_convergence}:} 1) Follows from the fact that \eqref{eq:merit} is a piecewise convex polynomial and therefore \cite[Theorem 4]{belousov2002frank} implies that the minimum is attained.

2) In Criterion~\ref{as:error_bound}, let $\partial h = \partial r^\star$ and therefore $S = [t^\star, \infty)$. Now, let $w \in \partial r^\star(t)$ for some fixed $t \in \R$ such that $|w| < \vartheta$. If $w > 0$, then Criterion~\ref{as:error_bound} holds vacuously since Proposition~\ref{prop:r_star} enforces $w \leq 0$. Now, assume $w \leq 0$. Following the definition of a subgradient for convex functions, $w \in \partial r^\star(t)$ if and only if $$r^\star(\tau) \geq r^\star(t) + w(\tau - t),$$
for all $\tau \in \R$. Now, let $\tau = \Pi_S(t) = t^\star$ and $M := \min_{x,s} r_0(x,s)$ and observe
$r^\star(t) \leq w(t - t^\star) + M.$ From Lemma~\ref{lem:quadratic-growth-rstar}, we have $\exists a > 0$ such that $a(t - t^\star)^2 \leq w(t - t^\star).$
Finally, we know from Proposition~\ref{prop:r_star} that $w \leq 0$ and $$|t - t^\star| \leq \frac{1}{a}|w|.$$
Therefore, Criterion~\ref{as:error_bound} is satisfied for any $\vartheta > 0$ and $\varphi = \frac{1}{a}$. 3) Follows as a direct result of Theorem~\ref{thm:exact_outer}.

\textit{Proof of Proposition~\ref{prop:inner}}: See \cite{ulbrich2011semismooth}, \cite{liao2020fbstab} and \cite{qi1993nonsmooth}.

\end{document}

%% file: plots/r_star.tex
\begin{tikzpicture}

\definecolor{darkgray176}{RGB}{176,176,176}
\definecolor{darkorange25512714}{RGB}{255,127,14}
\definecolor{lightgray204}{RGB}{204,204,204}
\definecolor{steelblue31119180}{RGB}{31,119,180}
\definecolor{forestgreen4416044}{RGB}{44,160,44}

\begin{axis}[
width = 0.9\linewidth,
height = 0.6\linewidth,
legend cell align={right},
legend style={fill opacity=0.5, draw opacity=1, text opacity=1, draw=lightgray204},
legend pos = north east,
tick align=outside,
tick pos=left,
title={\small Value Function $r^\star(\cdot)$},
x grid style={darkgray176},
xlabel={\small $t$},
xmajorgrids,
xmin=-0.7, xmax=1,
xtick style={color=black},
y grid style={darkgray176},
ylabel={\small $r^\star(t)$},
ymajorgrids,
ymin=-0.1, ymax=1,
ytick={0, 0.25, 0.5, 0.75, 1},
yticklabels={$0$, $0.25$, $0.5$, $0.75$, $1$},
ytick style={color=black}
]
\addplot [line width=0.5mm, steelblue31119180]
table {%
-1 1.33621353608209
-0.97979797979798 1.29409469620758
-0.95959595959596 1.25273912616321
-0.939393939393939 1.2121449554297
-0.919191919191919 1.17231023274846
-0.898989898989899 1.13323292232814
-0.878787878787879 1.09491089988768
-0.858585858585859 1.05734194853298
-0.838383838383838 1.02052375446543
-0.818181818181818 0.984453902520643
-0.797979797979798 0.949129871537488
-0.777777777777778 0.914549029558306
-0.757575757575758 0.880708628862733
-0.737373737373737 0.847605800839237
-0.717171717171717 0.815237550700275
-0.696969696969697 0.783600752049175
-0.676767676767677 0.752692141309211
-0.656565656565657 0.722508312028114
-0.636363636363636 0.693045709074209
-0.616161616161616 0.664300622743782
-0.595959595959596 0.636269182802843
-0.575757575757576 0.608947352490574
-0.555555555555556 0.582330922515795
-0.535353535353535 0.556415505082655
-0.515151515151515 0.531196527986271
-0.494949494949495 0.506669228824261
-0.474747474747475 0.482828649375171
-0.454545454545454 0.459669630200014
-0.434343434343434 0.437186805528315
-0.414141414141414 0.415374598494998
-0.393939393939394 0.394227216799123
-0.373737373737374 0.373738648859642
-0.353535353535353 0.353902660546861
-0.333333333333333 0.334712792570877
-0.313131313131313 0.316162358609868
-0.292929292929293 0.298244444261312
-0.272727272727273 0.280951906898004
-0.252525252525252 0.264277376507745
-0.232323232323232 0.24821325759078
-0.212121212121212 0.232751732182204
-0.191919191919192 0.217884764057612
-0.171717171717172 0.20360410416922
-0.151515151515151 0.189901297346557
-0.131313131313131 0.176767690280742
-0.111111111111111 0.164194440794602
-0.0909090909090908 0.15217252838268
-0.0707070707070706 0.14069276598597
-0.0505050505050504 0.129745812946512
-0.0303030303030303 0.119322189067182
-0.0101010101010101 0.109412289682879
0.0101010101010102 0.100006401631264
0.0303030303030305 0.0910947199949408
0.0505050505050506 0.082667365472977
0.0707070707070707 0.0747144022284101
0.0909090909090911 0.0672258560502748
0.111111111111111 0.0601917326639782
0.131313131313131 0.0536020360226555
0.151515151515152 0.0474467864144999
0.171717171717172 0.0417160382268127
0.191919191919192 0.0363998972164229
0.212121212121212 0.0314885371478294
0.232323232323232 0.0269722156744448
0.252525252525253 0.0228412893541989
0.272727272727273 0.0190862277079109
0.292929292929293 0.0156976262467483
0.313131313131313 0.0126662184132006
0.333333333333333 0.00998288639783093
0.353535353535354 0.00763867081119034
0.373737373737374 0.00562477920631928
0.393939393939394 0.0039325934619299
0.414141414141414 0.00255367604945307
0.434343434343434 0.00147977521850051
0.454545454545455 0.000702829144877365
0.474747474747475 0.000214969093078074
0.494949494949495 8.52165126068662e-06
0.515151515151515 4.83248392797712e-16
0.535353535353535 7.77198100509124e-17
0.555555555555556 3.01463804049018e-17
0.575757575757576 1.59191012102107e-17
0.595959595959596 9.85635848387328e-18
0.616161616161616 6.72423161915912e-18
0.636363636363636 4.90057982587459e-18
0.656565656565657 3.74222831243585e-18
0.676767676767677 2.962452832312e-18
0.696969696969697 2.41113837872189e-18
0.717171717171717 2.00685916455744e-18
0.737373737373737 1.70232517113469e-18
0.757575757575758 1.20047061006358e-18
0.777777777777778 9.85317352371341e-19
0.797979797979798 8.10156434373635e-19
0.818181818181818 6.66479277028713e-19
0.838383838383838 5.48115354610953e-19
0.858585858585859 4.50408174752558e-19
0.878787878787879 3.69722154281013e-19
0.898989898989899 3.03181393322912e-19
0.919191919191919 2.82106713616997e-20
0.939393939393939 2.82033142419693e-20
0.95959595959596 2.97453137156595e-20
0.97979797979798 3.12756565870651e-20
1 3.2787443070689e-20
};
\addlegendentry{Feasible}
\addplot [line width=0.5mm, darkorange25512714]
table {%
-1 1.58621353624191
-0.97979797979798 1.54409469636308
-0.95959595959596 1.50273912631483
-0.939393939393939 1.46214495557699
-0.919191919191919 1.42231023290934
-0.898989898989899 1.38323292247318
-0.878787878787879 1.34491090003432
-0.858585858585859 1.30734194868061
-0.838383838383838 1.27052375461342
-0.818181818181818 1.23445390267105
-0.797979797979798 1.19912987168945
-0.777777777777778 1.16454902970956
-0.757575757575758 1.13070862901432
-0.737373737373737 1.09760580099086
-0.717171717171717 1.0652375508511
-0.696969696969697 1.03360075219505
-0.676767676767677 1.00269214144592
-0.656565656565657 0.972508312140937
-0.636363636363636 0.943045709158199
-0.616161616161616 0.914300622797606
-0.595959595959596 0.886269183284034
-0.575757575757576 0.85894735258809
-0.555555555555556 0.832330922546323
-0.535353535353535 0.806415505118063
-0.515151515151515 0.781196528144586
-0.494949494949495 0.756669229389183
-0.474747474747475 0.732828649458848
-0.454545454545454 0.709669630320345
-0.434343434343434 0.687186805718127
-0.414141414141414 0.665374598694044
-0.393939393939394 0.644227216853568
-0.373737373737374 0.623738968395929
-0.353535353535353 0.603996084749135
-0.333333333333333 0.585069444538332
-0.313131313131313 0.566959047622877
-0.292929292929293 0.549664894120223
-0.272727272727273 0.533186984436445
-0.252525252525252 0.517525316952847
-0.232323232323232 0.502679892532389
-0.212121212121212 0.488650712305905
-0.191919191919192 0.475437774643714
-0.171717171717172 0.4630410803886
-0.151515151515151 0.451460629257596
-0.131313131313131 0.440696421395533
-0.111111111111111 0.430748458828689
-0.0909090909090908 0.421616736621554
-0.0707070707070706 0.413301258352858
-0.0505050505050504 0.405802023625985
-0.0303030303030303 0.399119031232615
-0.0101010101010101 0.393252282947497
0.0101010101010102 0.388201777907616
0.0303030303030305 0.383967516117779
0.0505050505050506 0.380549497595473
0.0707070707070707 0.377947722417509
0.0909090909090911 0.376162190094671
0.111111111111111 0.375192901236897
0.131313131313131 0.37500000000924
0.151515151515152 0.375000000023591
0.171717171717172 0.37500000000553
0.191919191919192 0.375000000010646
0.212121212121212 0.375000000025733
0.232323232323232 0.375000000039894
0.252525252525253 0.375000000053175
0.272727272727273 0.375000000065598
0.292929292929293 0.375000000077186
0.313131313131313 0.375000000087971
0.333333333333333 0.375000000097997
0.353535353535354 0.375000000107315
0.373737373737374 0.375000000115981
0.393939393939394 0.375000000124055
0.414141414141414 0.37500000013159
0.434343434343434 0.37500000013864
0.454545454545455 0.375000000145253
0.474747474747475 0.375000000151472
0.494949494949495 0.375000000157339
0.515151515151515 0.375000000162886
0.535353535353535 0.375000000168143
0.555555555555556 0.375000000173136
0.575757575757576 0.37500000017789
0.595959595959596 0.375000000182424
0.616161616161616 0.375000000186386
0.636363636363636 0.375000000036732
0.656565656565657 0.375000000037418
0.676767676767677 0.375000000038074
0.696969696969697 0.375000000038701
0.717171717171717 0.375000000039301
0.737373737373737 0.375000000039876
0.757575757575758 0.375000000040429
0.777777777777778 0.375000000040959
0.797979797979798 0.375000000041469
0.818181818181818 0.375000000041959
0.838383838383838 0.37500000004243
0.858585858585859 0.375000000042884
0.878787878787879 0.375000000043321
0.898989898989899 0.375000000043741
0.919191919191919 0.375000000044146
0.939393939393939 0.375000000044536
0.95959595959596 0.375000000044913
0.97979797979798 0.375000000045276
1 0.375000000045625
};
\addlegendentry{Infeasible}

\addplot [line width=0.5mm, forestgreen4416044]
table {%
-1.0 0.0
1.0 0.0
};
\addlegendentry{Unbounded}
\end{axis}

\end{tikzpicture}

%% file: plots/median_matrix.tex
\begin{tikzpicture}
  \begin{axis}[
    width=0.62\linewidth,
    height=0.62\linewidth,
    enlargelimits=false,
    title={\small Median Number of Newton Iterations},
    axis on top,
    ylabel={\small $\epsilon$--perturbation},
    xlabel={\normalsize $t_0$},
    colorbar,
    colormap/viridis,
    xtick={0,1,2,3,4},
    ytick={0,1,2,3,4,5},
    xticklabels={0.18,0.13,0.09,0.06,0.0},
    yticklabels={$10^{-6}$, $10^{-5}$, $10^{-4}$, $10^{-3}$, $10^{-2}$, \text{cold}},
  ]
    \addplot [
      matrix plot*,      
      point meta=explicit,
      mesh/cols=5        
    ]
    table [meta=z, row sep=\\] {
      x  y  z\\
    0 0 1.0 \\
1 0 3.0 \\
2 0 4.0 \\
3 0 4.0 \\
4 0 4.0 \\
0 1 2.0 \\
1 1 4.0 \\
2 1 5.0 \\
3 1 5.0 \\
4 1 5.0 \\
0 2 3.0 \\
1 2 5.0 \\
2 2 6.0 \\
3 2 6.0 \\
4 2 6.0 \\
0 3 5.0 \\
1 3 13.0 \\
2 3 14.0 \\
3 3 16.0 \\
4 3 15.0 \\
0 4 5.0 \\
1 4 13.0 \\
2 4 15.0 \\
3 4 17.0 \\
4 4 16.0 \\
0 5 8.0 \\
1 5 15.0 \\
2 5 16.0 \\
3 5 17.0 \\
4 5 16.0 \\
    };
  \end{axis}
\end{tikzpicture}

%% file: plots/feasible.tex

\begin{tikzpicture}
\definecolor{darkgray176}{RGB}{176,176,176}
\definecolor{darkorange25512714}{RGB}{255,127,14}
\definecolor{lightgray204}{RGB}{204,204,204}
\definecolor{steelblue31119180}{RGB}{31,119,180}

\begin{axis}[
width = 0.85\linewidth,
height = 0.65\linewidth,
legend cell align={right},
legend style={fill opacity=0.5, draw opacity=1, text opacity=1, draw=lightgray204, at={(0.97,0.72)}, anchor=north east},
tick align=outside,
tick pos=left,
title={\small Retrieving Feasible Solution from Suboptimal $t$},
x grid style={darkgray176},
xlabel={\small $t$},
xmajorgrids,
xmin=0.185, xmax=0.27,
xtick style={color=black},
y grid style={darkgray176},
ylabel={\small Median Newton Iterations},
ymajorgrids,
ymin=1, ymax=10,
ytick={2,3,4,5,6,7,8,9},
yticklabels={$2$,$3$,$4$,$5$,$6$,$7$,$8$,$9$},
ytick style={color=black}
]
    \addplot[color={rgb,1:red,0.0;green,0.6056;blue,0.9787}, name path={9}, draw opacity={1.0}, line width={1}, solid]
        table[row sep={\\}]
        {
            \\
            0.185538  9.0  \\
            0.18924876000000002  8.0  \\
            0.1930337352  7.0  \\
            0.196894409904  6.0  \\
            0.20083229810208  6.0  \\
            0.20484894406412163  5.0  \\
            0.20894592294540404  5.0  \\
            0.21312484140431215  5.0  \\
            0.2173873382323984  4.0  \\
            0.22173508499704636  4.0  \\
            0.22616978669698728  4.0  \\
            0.23069318243092704  4.0  \\
            0.23530704607954558  4.0  \\
            0.2400131870011365  4.0  \\
            0.24481345074115923  4.0  \\
            0.2497097197559824  4.0  \\
            0.2547039141511021  4.0  \\
            0.25979799243412416  4.0  \\
            0.2649939522828066  4.0  \\
            0.27029383132846274  4.0  \\
        }
        ;
    \addlegendentry {$\epsilon$ = 0.1}
    \addplot[color={rgb,1:red,0.8889;green,0.4356;blue,0.2781}, name path={10}, draw opacity={1.0}, line width={1}, solid]
        table[row sep={\\}]
        {
            \\
            0.185538  4.0  \\
            0.18924876000000002  4.0  \\
            0.1930337352  3.0  \\
            0.196894409904  4.0  \\
            0.20083229810208  3.0  \\
            0.20484894406412163  3.0  \\
            0.20894592294540404  4.0  \\
            0.21312484140431215  3.5  \\
            0.2173873382323984  4.0  \\
            0.22173508499704636  3.0  \\
            0.22616978669698728  3.0  \\
            0.23069318243092704  4.0  \\
            0.23530704607954558  3.5  \\
            0.2400131870011365  4.0  \\
            0.24481345074115923  3.0  \\
            0.2497097197559824  4.0  \\
            0.2547039141511021  4.0  \\
            0.25979799243412416  3.0  \\
            0.2649939522828066  4.0  \\
            0.27029383132846274  4.0  \\
        }
        ;
    \addlegendentry {$\epsilon$ = 0.01}
    \addplot[color={rgb,1:red,0.2422;green,0.6433;blue,0.3044}, name path={11}, draw opacity={1.0}, line width={1}, solid]
        table[row sep={\\}]
        {
            \\
            0.185538  2.0  \\
            0.18924876000000002  2.0  \\
            0.1930337352  2.0  \\
            0.196894409904  2.0  \\
            0.20083229810208  2.0  \\
            0.20484894406412163  2.0  \\
            0.20894592294540404  2.0  \\
            0.21312484140431215  2.0  \\
            0.2173873382323984  2.0  \\
            0.22173508499704636  2.0  \\
            0.22616978669698728  2.0  \\
            0.23069318243092704  2.0  \\
            0.23530704607954558  2.0  \\
            0.2400131870011365  2.0  \\
            0.24481345074115923  2.0  \\
            0.2497097197559824  2.0  \\
            0.2547039141511021  2.0  \\
            0.25979799243412416  2.0  \\
            0.2649939522828066  2.0  \\
            0.27029383132846274  2.0  \\
        }
        ;
    \addlegendentry {$\epsilon$ = 0.001}
    \addplot[color={rgb,1:red,0.7644;green,0.4441;blue,0.8243}, name path={12}, draw opacity={1.0}, line width={1}, solid]
        table[row sep={\\}]
        {
            \\
            0.185538  8.0  \\
            0.18924876000000002  8.0  \\
            0.1930337352  8.0  \\
            0.196894409904  8.0  \\
            0.20083229810208  8.0  \\
            0.20484894406412163  8.0  \\
            0.20894592294540404  8.0  \\
            0.21312484140431215  8.0  \\
            0.2173873382323984  8.0  \\
            0.22173508499704636  8.0  \\
            0.22616978669698728  8.0  \\
            0.23069318243092704  8.0  \\
            0.23530704607954558  8.0  \\
            0.2400131870011365  8.0  \\
            0.24481345074115923  8.0  \\
            0.2497097197559824  8.0  \\
            0.2547039141511021  8.0  \\
            0.25979799243412416  8.0  \\
            0.2649939522828066  8.0  \\
            0.27029383132846274  8.0  \\
        }
        ;
    \addlegendentry {cold start}
\end{axis}
\end{tikzpicture}